\documentclass{amsart}
\usepackage{amscd}
\usepackage{amsmath}
\usepackage{amsxtra}
\usepackage{amsfonts}
\usepackage{amssymb}

\newtheorem{theorem}{Theorem}[section]

\newtheorem{lemma}[theorem]{Lemma}
\newtheorem{proposition}[theorem]{Proposition}
\theoremstyle{definition}

\theoremstyle{remark}

\renewcommand{\theclaim}{\textup{\theclaim}}

\newtheorem*{acknowledgements}{Acknowledgements}

\numberwithin{equation}{section}

\def\openone

{\mathchoice

{\hbox{\upshape \small1\kern-3.3pt\normalsize1}}

{\hbox{\upshape \small1\kern-3.3pt\normalsize1}}

{\hbox{\upshape \tiny1\kern-2.3pt\SMALL1}}

{\hbox{\upshape \Tiny1\kern-2pt\tiny1}}}

\makeatletter

\newbox\ipbox

\newcommand{\diracb}[1]{\left\langle #1\mathrel{\mathchoice

{\setbox\ipbox=\hbox{$\displaystyle \left\langle\mathstrut 
#1\right.$}

\vrule height\ht\ipbox width0.25pt depth\dp\ipbox}

{\setbox\ipbox=\hbox{$\textstyle \left\langle\mathstrut 
#1\right.$}

\vrule height\ht\ipbox width0.25pt depth\dp\ipbox}

{\setbox\ipbox=\hbox{$\scriptstyle \left\langle\mathstrut 
#1\right.$}

\vrule height\ht\ipbox width0.25pt depth\dp\ipbox}

{\setbox\ipbox=\hbox{$\scriptscriptstyle \left\langle\mathstrut 
#1\right.$}

\vrule height\ht\ipbox width0.25pt depth\dp\ipbox}

}\right. }

\newcommand{\dirack}[1]{\left. \mathrel{\mathchoice

{\setbox\ipbox=\hbox{$\displaystyle \left.\mathstrut 
#1\right\rangle$}

\vrule height\ht\ipbox width0.25pt depth\dp\ipbox}

{\setbox\ipbox=\hbox{$\textstyle \left.\mathstrut 
#1\right\rangle$}

\vrule height\ht\ipbox width0.25pt depth\dp\ipbox}

{\setbox\ipbox=\hbox{$\scriptstyle \left.\mathstrut 
#1\right\rangle$}

\vrule height\ht\ipbox width0.25pt depth\dp\ipbox}

{\setbox\ipbox=\hbox{$\scriptscriptstyle \left.\mathstrut 
#1\right\rangle$}

\vrule height\ht\ipbox width0.25pt depth\dp\ipbox}

} #1\right\rangle}

\newcommand{\bz}{\mathbb{Z}}
\newcommand{\br}{\mathbb{R}}
\newcommand{\bc}{\mathbb{C}}
\newcommand{\bt}{\mathbb{T}}
\newcommand{\bn}{\mathbb{N}}

\newcommand{\xir}{X_\infty(r)}
\usepackage{graphicx}

\input cyracc.def

\begin{document}
\title[Disintegration of projective measures]{Disintegration of projective measures}
\author{Dorin Ervin Dutkay and Palle E.T. Jorgensen}
\address{Palle E.T. Jorgensen\\
Department of Mathematics\\
The University of Iowa\\
14 MacLean Hall\\
Iowa City, IA 52242-1419\\
U.S.A.\\}

\address{Dorin Ervin Dutkay\\
 Department of Mathematics\\
Hill Center-Busch Campus\\
Rutgers, The State University of New Jersey\\
110 Frelinghuysen Rd\\
Piscataway, NJ 08854-8019, USA }
 \email{Dorin Ervin Dutkay: ddutkay@math.rutgers.edu}\ 
\email{Palle E.T. Jorgensen: jorgen@math.uiowa.edu}
\thanks{Work supported in part by NSF grant DMS 0457491} \subjclass[2000]{42C40, 42A16, 42A65, 43A65, 46G15, 47D07, 60G18}
\keywords{Measures, projective limits, transfer operator, 
martingale, fixed-point}
\begin{abstract}
In this paper, we study a class of quasi-invariant measures on 
paths generated by discrete dynamical systems. Our main result 
characterizes the subfamily of these measures which admit a 
certain disintegration. This is a disintegration with respect to 
a random walk Markov process which is indexed by the starting 
point of the paths. Our applications include wavelet 
constructions on Julia sets of rational maps on the Riemann 
sphere.
\end{abstract}
 \maketitle 
\tableofcontents 
\section{\label{intro}Introduction}
We prove a general theorem about endomorphisms of measure spaces,
but it is motivated by our desire to use wavelet-type time/scale-self
similarity in operator theory, and more generally outside the traditional
context of wavelet theory; see e.g., \cite{BrJo02}. Hence our present setup is
chosen with these general applications in mind. What is typical is that the
starting point is geometric, and that a Hilbert space is not given at the
outset. The best choice of Hilbert space in turn is dictated by the
applications at hand. 
\par
Since wavelet decompositions use a special case of a general
disintegration process for measures via transition operators \cite{Rue89,Wal01},
we feel that it is of interest to offer a direct and self contained
statement and proof of a necessary and sufficient condition for
disintegration. 
\par
In its general form, the subject of disintegration of measures has its
roots in pioneering work by David Ruelle \cite{Rue89} on transfer operators. This
viewpoint was continued by Peter Walters \cite{Wal01} in ergodic theory; and in
operator algebras by a series of papers by Jean Renault and co-workers, on
the use of groupoids in representation theory and in dynamics, see \cite{Ren03,
KuRe03, AnRe01, Ren98, Ren80}. While our present setup and proofs are self
contained, in contrast the groupoid approach is realized via representations
of a certain non-abelian algebra of operators (a $C^*$-algebra), see also
\cite{BJO04}. Hence the other approaches use representations of non-abelian
algebras, which in turn entail more restrictive assumptions on the
endomorphism and on the transfer operator which is the starting point of the
analysis. Jean Renault's approach is via a groupoid $C^*$-algebra \cite{Ren80}. So
in this approach, the first step is from endomorphism to groupoid. 
\par
Our present approach is rather direct, starting with a given
endomorphism of a measure space. Yet our result Theorem 3.4 below belongs to
the same circle of ideas as that of a related disintgeration theorem of
Renault et al; see e.g., \cite[Theorems 2.8 and 3.5]{KuRe03}. However, turning to
our present setup, there are new elements, e.g., path space measures: We
begin with a construction of families of path space measures, and there is
no natural groupoid which is represented on our path space $X_\infty$ (see
section 3 below). So Renault's approach does not apply directly.
\par
Our starting point is an endomorphism in a measure space $X$. One might
add that it is then possible to get a corresponding representation of
Renault's groupoid built on $X$; or rather Renault's groupoid $C^*$-algebra.
With this, we note that our construction is then in fact a dilation of
Renault's representation. Having representations opens up for ways of doing
decomposition and disintegration; see also \cite{BJO04}. This viewpoint further
opens up some intriguing possibilities, and applications, which however we
leave to a future paper. 

\section{Endomorphisms and measures}

Let $X$ be a compact Hausdorff space, and let $r:X\rightarrow X$ 
be a given finite-to-one continuous endomorphism, mapping $X$ 
onto itself. \par
  As an example, $r = r(z)$ could be a rational mapping of the Riemann sphere. In 
that case, we take $X$ to be the corresponding Julia set; or $X$ 
could be the state-space of a subshift associated to a $0-1$ 
matrix. 
\par
Due to work of Brolin \cite{Bro} and Ruelle \cite{Rue89}, it is 
known that, for these examples, the dynamical system $(X,r)$ carries a unique maximal 
entropy, or minimal free energy measure $\mu$. Specifically, for 
each point $x$ in $X$, $\mu$ distributes the "mass" equally on 
the finite number of solutions $y$ to $r(y)=x$. The measure $\mu$ 
is obtained by iterating this over the paths of inverse images 
under $r$.
\par
  In earlier papers \cite{BrJo05}, \cite{DuJo04a},
\cite{DuJo04b}, we showed that this structure in fact admits a 
rich family of wavelet bases.  This wavelet construction is on
a Hilbert space which is different from, but analogous to, $L^2(\br)$, and
the more familiar case of multiwavelets on the real line.
\par
  These are the wavelets corresponding to scaling $x$ to $Nx$, by a fixed integer 
$N$, $N\geq 2$. In that case, $X=\bt$, the circle in $\bc$, and 
$r(z)=z^N$. So even in the "classical" case, there is a "unitary 
dilation" from $L^2(X)$ to $L^2(\br)$. In this construction, the 
Haar measure on $\bt$ "turns" into the Lebesgue measure on $\br$. 
In our general construction, we pass from $(X,r,\mu)$ to a 
Hilbert space built on a system $(X_\infty,\hat r,\hat\mu)$. Just 
as in the familiar case of $L^2(\br)$, we shall be concerned with 
the interplay between the initial space and the dilated one. Our 
consideration of quasi-invariant measures is motivated by wavelet 
filters in the more familiar case of $L^2(\br)$ wavelet analysis. 
See for example \cite{Dau92} and \cite{BrJo02}. We stress that it 
is the dilated Hilbert space which carries wavelet bases. These 
bases are in turn parametrized by a family of functions $V$ on 
$X$, or on $\bt$ in the classical case.
\par
Our work so far, on examples shows that this viewpoint holds 
promise for a deeper understanding the harmonic analysis of Julia 
sets, and of other iteration systems. 
\par
  In these cases, the analogue of $L^2(\br)$ involves quasi-invariant 
measures $\hat\mu$ on a space $X_\infty$, built from $X$ in a way 
that follows closely the analogue of passing from $\bt$ to $\br$ 
in wavelet theory; see \cite{DuJo04a} for details. But the 
induced measure $\hat\mu$ on $X_\infty$ is now a path-space 
measure. The analogue of the unitary action of translations by 
the integers $\bz$ on $L^2(\br)$, or equivalently multiplication 
in the Fourier dual, will be described in general below. For the 
Julia examples, it is the operator of multiplication by $z$ on 
$L^2(X,\mu)$, i.e., it is a normal operator. We get a 
corresponding covariant system on $X_\infty$, where 
multiplication by $f(z)$ is unitarily equivalent to 
multiplication by the composition $f(r(z))$. But this is now 
dilated onto the Hilbert space $L^2(X_\infty)$, and our dilated 
Hilbert space is defined from our path-space measures $\hat\mu$.
\par
  Hence all the issues from wavelet analysis of the mid nineteen eighties, especially the pyramid algorithms, have analogues 
in this wider context. In other work, we use this to construct multiresolutions on 
the Julia sets, and on other iteration systems, and our $r$-scale 
multiresolutions appear to reveal interesting spectral theory for 
Julia sets. And it seems further to offer a promising interplay 
between geometry and spectrum. None of this is yet well 
understood.
\par
  The purpose of this paper is to give an interconnection between the measures $\mu_0$ on $X$, and 
a family of $V$-quasi-invariant measures $\hat\mu$ on $X_\infty$. 
We characterize those $V$-quasi-invariant measures $\hat\mu$ on 
$X_\infty$ which admit a path-space decomposition in terms of the 
initial measure $\mu_0$ at the start of the paths. This is useful because it gives the $V$-quasi-invariant measure a more concrete form, see \cite{DuJo04a}.
\par
We recall the definitions: if $\mathfrak{B}$ is a sigma-algebra 
on a set $X$ and $r:X\rightarrow X$ is a finite-to-one, onto and 
measurable map, we denote by
$$\xir:=\{(x_0,x_1,...)\in\prod_{n\in\bn_0}X\,|\,r(x_{n+1})=x_n,\mbox{ 
for all }n\in\bn_0\}.$$ We denote the projections by
$\theta_n:\xir\ni(x_0,x_1,...)\mapsto x_n\in X.$ The union of the 
pull-backs $\theta_n^{-1}(\mathfrak{B})$ generates a 
sigma-algebra $\mathfrak{B}_\infty$. The map $r$ extends to a 
measurable bijection $\hat r:\xir\rightarrow\xir$ defined by
$$\hat r(x_0,x_1,...)=(r(x_0),x_0,x_1,...).$$
\par
Let $V:X\rightarrow[0,\infty)$ be a measurable, bounded function. 
We say that a measure $\mu_0$ on $(X,\mathfrak{B})$ has the {\it 
fixed-point property} if
\begin{equation}\label{eq2_1_1}
\int_XV\,f\circ r\,d\mu_0=\int_Xf\,d\mu_0,\quad(f\in 
L^\infty(X,\mu_0).\end{equation}
 We say that a measure $\hat\mu$ on 
$(\xir,\mathfrak{B}_\infty)$ is {\it $V$-quasi-invariant} if 
$$d(\hat\mu\circ\hat r)=V\circ\theta_0\,d\hat\mu.$$
We recall the result from \cite{DuJo04b}:
\begin{theorem}\label{th1_1}
There exists a one-to-one correspondence between measures $\mu_0$ 
on $X$ with the fixed-point property and $V$-quasi-invariant 
measures $\hat\mu$ on $\xir$, given by
$$\mu_0=\hat\mu\circ\theta_0^{-1}.$$
\end{theorem}

Our main result Theorem \ref{th2_4} below applies to a class of measures
on some extended spaces $X_\infty$, generalized solenoids. Our setup is as
follows: Start with $(X, r, V, \mu_0)$ where $X$ is a measure space with an
endomorphism $r$, $V$ is a weight function as specified, and $\mu_0$ is a
$V$-strongly invariant measure. In wavelet applications, see \cite{DuJo04b}, $V$ may
have the form $V = |m_0|^2$ where $m_0$ is a wavelet low-pass filter. In the
dyadic case, $X = \bt$ (the circle or the one-torus), and  $r: z \mapsto z^2$. Under
suitable conditions on the so called filter function $m_0$, we show in
\cite{DuJo04b} that $(X_\infty, \mu_\infty)$ is measure theoretically equivalent
to the real line $(\br,\mbox{ Lebesgue measure })$; and so, in this case, the dilated
Hilbert space is merely the familiar $L^2(\br)$ of classical wavelet theory
\cite{Dau92}. In \cite{DuJo04a}, we show that the same idea applies to the case when
the system $(X, r, V, \mu_0)$ is as follows: $X:= X(r)$ is the Julia set of some
fixed rational mapping $r(z)$ of one complex variable. (We denote by $r$ also
the restriction to $X$ of the rational mapping $r(z)$.)  Using a theorem of
Brolin \cite{Bro} about measures on Julia sets, we further describe in \cite{DuJo04a}
a class of weight functions $V$ and measures $\mu_0$ on $X$ which satisfy the
$V$-fixed point property considered in Theorem \ref{th2_1}. We also describe in
\cite{DuJo04a} how our dilation construction then provides wavelet bases in the
dilation Hilbert space $L^2(X_\infty, \mu_\infty)$. However, in our work in
\cite{DuJo04a,DuJo04b} some cases are left open. This limits our approach to
dilation/wavelet constructions for dynamical systems. The reason is that, at
the time we lacked a general disintegration formula for the relevant measure
spaces. It is the purpose of the present paper to remedy this, and to prove
such a theorem, specifically Theorem \ref{th2_4} below. We now turn to the
path-space measures which are used in our theorem.

\par
The significance of our path-space measures $\{P_x\, |\, x\in X\}$ is that
they support our disintegration from Theorem \ref{th2_4} 
below. This is further
elaborated in \cite{DuJo04a,DuJo04b}. In the wavelet case, starting with $V = |m_0|^2$,
let $P_x$ be the path space measures from Proposition \ref{prop2_3}.  We gave conditions
in \cite{DuJo04a} which imply that for a. e. $x$, $P_x$ is supported on the union $U$
of certain finite orbits for $r$. We outline how this union $U$ is naturally in
a bijective correspondence with a copy of the integers $\bz$. This turns out to
be consistent with the natural way the group $\bz$ acts by translation in
wavelet theory of the Hilbert space $L^2(\br)$.
\par
Starting with $V$ and the corresponding $P_x$, we then get an
important function $h(x):= P_x(\bz)$. It turns out to be a fundamental
eigenfunction of wavelet theory. In fact the function $h$ is known \cite{BrJo02} to
be a fundamental minimal eigenfunction for the wavelet transfer operator $R_V$
(also called a Ruelle-Perron-Frobenius operator), studied first by W. Lawton
in \cite{Law90}. The function $h$ determines a number of important properties of a
wavelet basis; see e.g., \cite[Chapter 6]{Dau92} and \cite{BrJo02} for additional
details on this point.  We further show in \cite{DuJo04b} that $R_V h = h$ holds,
and further that $h$ is an infinite-dimensional generalization of the familiar
Perron-Frobenius right-eigenvector. This entails a suitable renormalization
so that the Perron-Frobenius eigenvalue is $1$. The measures $\mu_0$ from Theorem
\ref{th2_1} satisfying our $V$-fixed point property may then be viewed as
generalizations of left-Perron-Frobenius eigenvectors. See also \cite{BrJo02} and
\cite{Rue89} for details on infinite-dimensional Perron-Frobenius theory.

\section{\label{disi}Disintegration}

\par
Let $X$ be a non-empty set, $\mathfrak{B}$ a sigma-algebra of 
subsets of $X$ and $r:X\rightarrow X$, an onto, finite-to-one, and
$\mathfrak{B}$-measurable map. 
\par
We will assume in addition that we can label measurably the 
branches of the inverse of $r$. By this, we mean that the 
following conditions are satisfied:
\begin{equation}\label{eq2_1}
\mbox{The map}\quad \mathfrak{c}:X\ni 
x\mapsto\#r^{-1}(x)<\infty\mbox{ is measurable }.
\end{equation}
We denote by
$$A_i:=\{x\in X\,|\, \mathfrak{c}(x)=\#r^{-1}(x)\geq 
i\},\quad(i\in\bn).$$ Equation (\ref{eq2_1}) implies that the 
sets $A_i$ are measurable. Also they form a decreasing sequence 
and, since the map is finite to one
$$X=\cup_{i=1}^\infty(A_{i+1}\setminus A_i).$$
Then, we assume that there exist measurable maps 
$\tau_i:A_i\rightarrow X$, $i\in\{1,2...\}$ such that 
\begin{equation}\label{eq2_2}
r^{-1}(x)=\{\tau_1(x),....,\tau_{\mathfrak{c}(x)}(x)\},\quad 
(x\in X),
\end{equation}
\begin{equation}\label{eq2_3}
\tau_i(A_i)\in\mathfrak{B},\mbox{ for all }i\in\{1,2...,\}.
\end{equation}
Thus $\tau_1(x),...,\tau_{\mathfrak{c}(x)}(x)$ is a list without 
repetitions of the "roots" of $x$, $r^{-1}(x)$.
\par
>From equation (\ref{eq2_2}) we obtain also that 
\begin{equation}\label{eq2_4}
\tau_i(A_i)\cap\tau_j(A_j)=\emptyset,\mbox{ if }i\neq j,
\end{equation}
and
\begin{equation}\label{eq2_5}
\cup_{i=1}^\infty\tau_i(A_i)=X.
\end{equation}
\par In the sequel, we will use the following notation: for a 
function $f:X\rightarrow\bc$, we denote by 
$$f\circ\tau_i(x):=\left\{\begin{array}{ccc}
f(\tau_i(x)),&\mbox{if}&x\in A_i,\\
0,&\mbox{if}&x\in X\setminus A_i.\end{array}\right.
$$
\begin{theorem}\label{th2_1}
Let $(X,\mathfrak{B})$ and $r:X\rightarrow X$ be as above. Let 
$V:X\rightarrow [0,\infty)$ be a bounded 
$\mathfrak{B}$-measurable map. For a measure $\mu_0$ on 
$(X,\mathfrak{B})$, the following affirmations are equivalent 
\begin{enumerate}
\item The measure $\mu_0$ has the fixed-point property (\ref{eq2_2_1}).
\item There exists a non-negative, $\mathfrak{B}$-measurable map 
$\Delta$ (depending on $V$ and $\mu_0$) such that 
\begin{equation}\label{eq2_1_2}
\sum_{r(y)=x}\Delta(y)=1,\mbox{ for }\mu_0-\mbox{a.e. }x\in X,
\end{equation}
and
\begin{equation}\label{eq2_1_3}
\int_XV\,f\,d\mu_0=\int_X\sum_{r(y)=x}\Delta(y)f(y)\,d\mu_0(x),\mbox{ 
for all }f\in L^\infty(X,\mu_0).
\end{equation}
\end{enumerate}
Moreover, when the affirmations are true, $\Delta$ is unique up 
to $V\,d\mu_0$-measure zero.
\end{theorem}

\begin{proof}
We check first the implication (ii)$\Rightarrow$(i). We have, for 
$f\in L^\infty(X,\mu_0)$:
\begin{align*}
\int_XV\,f\circ 
r\,d\mu_0&=\int_X\sum_{r(y)=x}\Delta(y)f(r(y))\,d\mu_0(x)\\
&=\int_Xf(x)\sum_{r(y)=x}\Delta(y)\,d\mu_0(x)=\int_X f\,d\mu_0,
\end{align*}
and (\ref{eq2_1_1}) follows.
\par
For the implication (i)$\Rightarrow$(ii) we will need a lemma:
\begin{lemma}\label{lem2_2}
Let $\mu_0$ be a measure on $(X,\mathfrak{B})$ that satisfies 
(\ref{eq2_1_1}). Then, for each $f\in L^1(X,\mu_0)$ there exists 
a function which we denote by $R_{\mu_0}(V\,f)\in L^1(X,\mu_0)$ 
such that
\begin{equation}\label{eq2_2_1}
\int_XV\,f\,g\circ 
r\,d\mu_0=\int_XR_{\mu_0}(V\,f)\,g\,d\mu_0,\quad(g\in 
L^\infty(X,\mu_0)).
\end{equation}
The operator $f\mapsto R_{\mu_0}(V\,f)$ has the following 
properties:
\begin{equation}\label{eq2_2_2}
R_{\mu_0}(V)(=R_{\mu_0}(V\cdot\mathbf{1}))=\mathbf{1},\mbox{ 
pointwise }\mu_0-\mbox{a.e.}
\end{equation}
\begin{equation}\label{eq2_2_3}
R_{\mu_0}(V\,f\,h\circ r)=R_{\mu_0}(V\,f)\,h,\quad(f,h\in 
L^\infty(X,\mu_0)).
\end{equation}
\begin{equation}\label{eq2_2_4}
R_{\mu_0}(Vf)\geq 0\mbox{ if }f\geq 0.
\end{equation}

\end{lemma}
\begin{proof}{\it of lemma \ref{lem2_2}}.
The positive functional 
$$g\mapsto \int_XV\,f\,g\circ r\,d\mu_0,\quad(g\in 
L^\infty(X,\mu_0)),$$ defines a measure on $(X,\mathfrak{B})$ 
which is absolutely continuous with respect to $\mu_0$; indeed, 
if $E\in\mathfrak{B}$ has $\mu_0(E)=0$ then, by (\ref{eq2_1_1})
$$\int_X V\,\chi_E\circ r\,d\mu_0=\int_X\chi_E\,d\mu_0=0,$$
therefore $r^{-1}(E)$ has $V\,d\mu_0$-measure zero so 
$$\int_XV\,f\,\chi_E\circ r\,d\mu_0=0.$$
Denoting by $R_{\mu_0}(V\,f)$ the Radon-Nikodym derivative, the 
relation (\ref{eq2_2_1}) is obtained.
\par
Taking $f=\mathbf{1}$, we obtain, for all $g\in 
L^\infty(X,\mu_0)$, 
$$\int_X R_{\mu_0}(V)\,g\,d\mu_0=\int_XV\,g\circ r\,d\mu_0=\int_X 
g\,d\mu_0,$$ where we used (\ref{eq2_1_1}) in the last step. Thus 
(\ref{eq2_2_2}) follows.
\par
Next, for all $g\in L^\infty(X,\mu_0)$ we have
$$\int_XR_{\mu_0}(V\,f\,h\circ r)\,g\,d\mu_0=\int_XV\,f\,h\circ r\,g\circ 
r\,d\mu_0=\int_XR_{\mu_0}(Vf)\,h\,g\,d\mu_0,$$ and this implies 
(\ref{eq2_2_3}).
\par
Take now $f\geq0$, $f\in L^1(X,\mu_0)$. For all $g\in 
L^\infty(X,\mu_0)$ with $g\geq0$, we have
$$\int_XR_{\mu_0}(Vf)\,g\,d\mu_0=\int_X V\,f\,g\circ 
r\,d\mu_0\geq0,$$ and (\ref{eq2_2_4}) follows.
\end{proof}
We return to the proof of the theorem. We prove that there exist 
some non-negative, measurable function $\Delta$ on $X$, such that 
\begin{equation}\label{eq2_1_4}
(R_{\mu_0}(V\,f))(x)=\sum_{i=1}^\infty\Delta(\tau_i(x))f(\tau_i(x)),\mbox{ 
for }\mu_0-\mbox{a.e. }x\in X,\quad(f\in L^1(X,\mu_0)).
\end{equation}
Note that every $f\in L^1(X,\mu_0)$ can be written as
\begin{equation}\label{eq2_1_5}
f=\sum_{i=1}^N(f\circ\tau_i)\circ r\,\chi_{\tau_i(A_i)}.
\end{equation}
Indeed, if $x\in X$, then, by (\ref{eq2_5}), $x\in\tau_i(A_i)$ 
for some $i\in\{1,2,...\}$, and $$x\in 
r^{-1}(r(x))=\{\tau_1(r(x)),...,\tau_{\mathfrak{c}(r(x))}(r(x))\}.$$ 
Therefore $x=\tau_i(r(x))$ and (\ref{eq2_1_5}) follows. Note also 
that, for all $x\in X$, the sum in (\ref{eq2_1_5}) is actually 
finite. It has at most $\mathfrak{c}(x)$ terms.
\par
The function $f\circ\tau_i$ is measurable 
because $\tau_i$ is.
\par
Applying $R_{\mu_0}$ and use (\ref{eq2_2_3}) we obtain that 
\begin{equation}\label{eq2_1_6}
R_{\mu_0}(Vf)=\sum_{i=1}^\infty 
R_{\mu_0}(V\,\chi_{\tau_i(A_i)})\,f\circ\tau_i.
\end{equation}
\par
A little argument is needed for this, to deal with the infinite 
sum. If $S_n$ are the partial sums 
$S_n:=\sum_{i=1}^n\chi_{\tau_i(A_i)}\,f\circ\tau_i\circ r$, then 
they converge pointwise and dominated to the limit $f$. The 
definition of $R_{\mu_0}$ shoes that $R_{\mu_0}(V\,S_n)$ 
converges weakly to $R_{\mu_0}(V\,f)$. However, the sum in 
$R_{\mu_0}(V\,f)$ is pointwise finite, because any $x\in X$ 
belongs to only finitely many $A_i$'s, so $f\circ\tau_i(x)=0$, 
except for finitely many $i$'s. Hence (\ref{eq2_1_6}) holds.
 \par Define then 
\begin{equation}\label{eq2_1_7}
\Delta(\tau_i(x))=(R_{\mu_0}(V\,\chi_{\tau_i(A_i)}))(x),\quad(x\in 
A_i)
\end{equation}
Since $(\tau_i(X))_{i=1,\infty}$ form a partition of $X$, 
equation (\ref{eq2_1_7}) defines the function $\Delta$ on $X$. 
\par
We claim that $\Delta$ is measurable. Take $E$ a Borel subset of 
$\br$. Then
\begin{align*}
\Delta^{-1}(E)&=\cup_{i=1}^\infty(\Delta^{-1}(E)\cap\tau_i(A_i))\\
&=\cup_{i=1}^\infty\left(\tau_i\left((\Delta\circ\tau_i)_{|_{A_i}}^{-1}(E)\right)\right)\\
&=\cup_{i=1}^\infty\left(\tau_i\left((R_{\mu_0}(V\chi_{\tau_i(A_i)}))^{-1}(E)\right)\right)\\
&=\cup_{i=1}^\infty\left(r^{-1}\left((R_{\mu_0}(V\chi_{\tau_i(A_i)}))^{-1}(E)\right)\cap 
\tau_i(A_i)\right),\end{align*} and since the functions 
$R_{\mu_0}(V\chi_{\tau_i(A_i)})$ are measurable, it follows that 
$\Delta$ is too.
\par
Now the equations (\ref{eq2_1_6}) and (\ref{eq2_1_7}) imply 
(\ref{eq2_1_4}). Plugging (\ref{eq2_1_4}) into (\ref{eq2_2_1}) we 
obtain (\ref{eq2_1_3}).
\par
$$\sum_{i=1}^\infty\Delta\circ\tau_i=R_{\mu_0}(V\,\mathbf{1})=1,$$
where we used (\ref{eq2_2_2}) in the last step. This proves 
(\ref{eq2_1_2}).
\par
To prove the uniqueness, let $\Delta'$ be another such function. 
Considering in (\ref{eq2_1_3}) functions $f$ supported on 
$\chi_{\tau_i(X)}$, we obtain that 
$$\Delta(\tau_i(x))=\Delta'(\tau_i(x)),\,\mbox{ for 
}\mu_0-\mbox{a.e. }x\in A_i.$$ But \begin{align*} \{x\in 
X\,|\,\Delta(x)\neq\Delta'(x)\}&=\cup_{i=1}^\infty\{x\in\tau_i(A_i)\,|\,\Delta(\tau_i(r(x)))\neq\Delta'(\tau_i(r(x)))\}\\
&=\cup_{i=1}^\infty\left(\tau_i(A_i)\cap r^{-1}\left(\{y\in X
,|\,\Delta(\tau_i(y))\neq\Delta'(\tau_i(y))\}\right)\right)
\end{align*}
But the argument in the beginning of the proof shows that if a 
set $E$ has $\mu_0(E)=0$ then $r^{-1}(E)$ has $V\,d\mu_0$-measure 
$0$. Thus each term of the union has $V\,d\mu_0$-measure $0$ and 
the theorem is proved.
\end{proof}

\subsection{\label{path}Path measures}
Each point in $(x_0,x_1,...)\in\xir$ is determined by a point 
$x_0\in X$ and a choice of a root $x_1\in r^{-1}(x_0)$, a root 
$x_2\in r^{-1}(x_1)\subset r^{-2}(x_0)$ and so on. \par If the 
point $x_0$ is fixed, then to specify a point in $\xir$ is to 
specify the choice of roots, or the choice of a path 
$$x_0\stackrel{r}{\longleftarrow}x_1\stackrel{r}{\longleftarrow}...\stackrel{r}{\longleftarrow} 
x_n\stackrel{r}{\longleftarrow}x_{n+1}\stackrel{r}{\longleftarrow}...$$ We will denote the 
set of these paths by 
$$\Omega_{x_0}:=\{(x_1,x_2,...)\,|\,(x_0,x_1,x_2,...)\in\xir\}.$$
The set can be regarded as a subset of $\xir$, namely the subset 
of points in $\xir$ that have the first coordinate equal to 
$x_0$. It can be also be regarded as the projective limit of the 
following diagram
$$
r^{-1}(x_0)\stackrel{r}{\longleftarrow} 
r^{-2}(x_0)\stackrel{r}{\longleftarrow}...\stackrel{r}{\longleftarrow}r^{-n}(x_0)\stackrel{r}{\longleftarrow} 
r^{-(n+1)}(x_0)\stackrel{r}{\longleftarrow}...
$$
\par
We can construct a Radon measure $P_{x_0}$ on $\Omega_{x_0}$ 
using $\Delta$ to assign probabilities to the choices of roots.
\begin{proposition}\label{prop2_3}
Let $(X,\mathfrak{B})$, $r:X\rightarrow X$ and be as above, and 
let $D:X\rightarrow[0,\infty)$ be a measurable function with the 
property that
\begin{equation}\label{eq2_3_1}
\sum_{r(y)=x}D(y)=1.
\end{equation}
Denote by 
\begin{equation}\label{eq2_3_2}
D^{(n)}:=D\cdot D\circ r\cdot...\cdot D\circ 
r^{n-1},\quad(n\in\bn),\quad D^{(0)}:=1.
\end{equation} Then for each $x_0\in X$, there exists a Radon probability measure $P_{x_0}$ 
on $\Omega_{x_0}$ such that, if $f$ is a bounded measurable 
function on $\Omega_{x_0}$ which depends only on the first $n$ 
coordinates $x_1,...,x_n$, then
\begin{equation}\label{eq2_3_3}
\int_{\Omega_{x_0}}
f(\omega)\,dP_{x_0}(\omega)=\sum_{r^n(x_n)=x_0}D^{(n)}(x_n)f(x_1,...,x_n).
\end{equation}
\end{proposition}
\begin{proof}
To check the consistency, take a measurable function $f$ which 
depends only on the first $n$ coordinates $x_1,...,x_n$. We 
regard it as a function which depends on the first $n+1$ 
coordinates and check that the two formulas in (\ref{eq2_3_3}) 
(one for $n$, one for $n+1$) give the same result:
\begin{align*}
\int_{\Omega_{x_0}}f(\omega)\,dP_{x_0}(\omega)&=\sum_{r^{n+1}(x_{n+1})=x_0}D^{(n+1)}(x_{n+1})f(x_1,...,x_{n+1})\\
&=\sum_{r^n(x_n)=x_0}f(x_1,...,x_n)D^{(n)}(x_n)\sum_{r(x_{n+1})=x_n} D(x_{n+1})\\
&=\sum_{r^n(x_n)=x_0}D^{(n)}(x_n)f(x_1,...,x_n),
\end{align*}
so $P_{x_0}$ is well defined. Since the functions which depend 
only on finitely many coordinates form a dense subalgebra of 
bounded functions (according to the Stone-Weierstrass theorem), 
we can use Riesz' theorem to obtain the probability measure 
$P_{x_0}$ on $\Omega_{x_0}$. The application of Stone-Weierstrass is justified because when the
point x is fixed, we may use Tychonoff-compactness of the usual infinite
product $\Omega:= \bz_N \times\bz_N\times...$. But since we may have a variable
number of preimages, we must first construct a dense subalgebra of bounded
functions, and then the measures $(P_x)$ by using the structure of the measure
space $X_{\infty}(r)$ defined above.
\end{proof}

\subsection{Main theorem}
\par
The next theorem shows that the $V$-quasi-invariant measure 
$\hat\mu$ on $\xir$ associated to $\mu_0$ can be disintegrated 
through the measures $P_{x_0}$ associated to the function 
$\Delta$.

\begin{theorem}\label{th2_4}
Let $(X,\mathfrak{B})$, $r:X\rightarrow X$ and 
$V:X\rightarrow[0,\infty)$ be as above. Let $\mu_0$ be a measure 
on $(X,\mathfrak{B})$ with the fixed point property 
(\ref{eq2_1_1}). Let $\Delta:X\rightarrow[0,1]$ be the function 
associated to $V$ and $\mu_0$ as in theorem \ref{th2_1}, and let 
$\hat\mu$ be the unique $V$-quasi-invariant measure on $\xir$ 
with 
$$\hat\mu\circ\theta_0^{-1}=\mu_0,$$
as in theorem \ref{th1_1}. For $\Delta$, we define the measures 
$P_{x_0}$ as in proposition \ref{prop2_3}. Then, for all bounded 
measurable functions $f$ on $\xir$,
\begin{equation}\label{eq2_4_1}
\int_{\xir}f\,d\hat\mu=\int_X\int_{\Omega _{x_0}}
f(x_0,\omega)\,dP_{x_0}(\omega)\,d\mu_0(x_0).
\end{equation}
\end{theorem}
\begin{proof}
By density, it is enough to prove the theorem for functions $f$ 
on $\xir$ which depend only on finitely many coordinates. So, 
assume $f=g\circ\theta_n$, with $g$ bounded and $\mathfrak{B}$ 
measurable. Then $g\circ\theta_n$ will depend only on $x_0$ and 
$x_1,...,x_{n}$,
$$g\circ\theta_n(x_0,x_1,...,x_n)=g(x_n).$$
\par
>From \cite[Theorem 3.3]{DuJo04b}, we recall that, if 
$\mu_n:=\hat\mu\circ\theta_n^{-1}$ then $d\mu_n=V^{(n)}\,d\mu_0$. 
Thus
$$\int_{\xir}g\circ\theta_n\,d\hat\mu=\int_Xg\,d\mu_n=\int_XV^{(n)}\,g\,d\mu_0.$$
Applying lemma \ref{lem2_2}, we obtain by induction that
$$\int_XV^{(n)}\,g\,d\mu_0=\int_XV\,g\,V^{(n-1)}\circ 
r\,d\mu_0=\int_XR_{\mu_0}(V\,g)V^{(n-1)}\,d\mu_0$$$$=...=\int_XR_{\mu_0}^n(V\,g)\,d\mu_0.$$ 
Also, using induction on (\ref{eq2_1_4}), we obtain that
$$(R_{\mu_0}^n(V\,g))(x_0)=\sum_{r^n(x_n)=x_0}\Delta^{(n)}(x_n)g(x_n)$$$$=\int_{\Omega_{x_0}}(g\circ\theta_n)(x_0,\omega)\,dP_{x_0}(\omega).$$
Combining these equalities we obtain (\ref{eq2_4_1}).
\end{proof}
\begin{acknowledgements}
We are pleased to thank an anonymous referee for calling our attention
to the connection to groupoids, and in particular to the work of Jean
Renault et al, stressing that this is an important alternative approach to
disintegration based on transfer operators. In addition, both authors
acknowledge encouragement and helpful discussions with David Larson; and
partial support from the National Science Foundation, DMS 0139473, and DMS
0457581.
\end{acknowledgements}

\end{document}